\documentclass{article}
\usepackage{amssymb}
\usepackage{hyperref}
\newcommand{\MR}[1]{\href{http://www.ams.org/mathscinet-getitem?mr=#1}{MR #1}}
\newcommand{\arXiv}[1]{\href{http://front.math.ucdavis.edu/#1}{arXiv:#1}}
\newcommand{\eq}{\begin{equation}}
\newcommand{\en}{\end{equation}}

\newcommand{\re}[1]{\mbox{(\ref{#1})}}

\newcommand{\prob}{\mathbb P}

\newcommand{\W}{W}
\newcommand{\V}{V}

\def\endpf{\hfill $\Box$ \vskip0.5cm}

\newtheorem{theorem}{\large Theorem}

\newtheorem{definition}[theorem]{\large Definition}
\newtheorem{corollary}[theorem]{\large Corollary}
\newtheorem{lemma}[theorem]{\large Lemma}

\begin{document}
\title{Exchangeable Gibbs partitions and Stirling triangles
\thanks{Research supported in part by N.S.F. Grant DMS-0405779}
}
\author{Alexander Gnedin\thanks{Utrecht University; e-mail gnedin@math.uu.nl}
\hspace{.2cm}
and 
\hspace{.2cm}
Jim Pitman\thanks{University of California, Berkeley; e-mail pitman@stat.Berkeley.EDU} 
\\
\\
\\
\\
}
\date{
\today
\\
}
\maketitle

\centerline{\bf Abstract}
\noindent 
For two collections of 
nonnegative and suitably normalised
weights $\W=(\W_j)$ and $\V=(\V_{n,k})$, 
a probability distribution on the set of partitions of the set $\{1,\ldots,n\}$ 
is defined by assigning to a generic partition $\{A_j, j\leq k\}$ the probability 
$\V_{n,k}\,\W_{|A_1|}\cdots \W_{|A_k|}$, where $|A_j|$ is the number of elements of $A_j$.
We impose constraints on the weights by assuming that 
the resulting random partitions $\Pi_n$ of $[n]$ are consistent as $n$ varies,
meaning that they define an exchangeable partition of the set of all natural numbers.
This implies that the 
weights $\W$ must be of a very special form depending on a single parameter
$\alpha\in [-\infty,1]$.
The case $\alpha=1$ is trivial, 
and for each value of $\alpha\neq 1$ the set of possible $\V$-weights is an infinite-dimensional simplex.
We identify the extreme points of the simplex by
solving the boundary problem for
a generalised Stirling triangle. In particular, we show that 
the boundary is discrete for $-\infty\leq\alpha<0$ and continuous for $0\leq\alpha<1$.
For $\alpha\leq 0$ the extremes correspond to the members of the Ewens-Pitman family of 
random partitions
indexed by $(\alpha,\theta)$, while for 
$0<\alpha<1$ the extremes are obtained by conditioning an $(\alpha,\theta)$-partition
on the asymptotics of the number of blocks of $\Pi_n$ as $n$ tends to infinity.

\vskip0.5cm

\section{Introduction}
\label{intro}

By a random partition of the set of natural numbers ${\mathbb N}$ we mean a consistent sequence $\Pi=(\Pi_n)$
of random partitions of finite sets
$[n]:=\{1,\ldots, n\}$. 
For each $n$ the range of the random variable $\Pi_n$ 
is the set of all partitions of $[n]$
into some number of disjoint nonempty blocks, 
and the consistency means that
$\Pi_{n}$ is obtained from $\Pi_{n+1}$ by discarding the element $n+1$.
A random partition $\Pi$ is {\em exchangeable} if for each $n$ the probability distribution of $\Pi_n$
is invariant under all permutations of $[n]$.

\par Let $\{A_j,\,1 \le j\leq k\}$ denote a generic partition of the set $[n]$, 
and let the $A_j$ be indexed by $[k]$ in order of their least elements.
Exchangeability of $\Pi$
means that 
$$
{\mathbb P}(\Pi_n=\{A_1,\ldots, A_k\})= p(|A_1|, \ldots, |A_k|)
$$
for some nonnegative function 
$$
p(\lambda) := p (\lambda_1,\ldots,\lambda_k)
$$
of {\em compositions} $\lambda=(\lambda_1,\ldots,\lambda_k)$ of $n$, such that
$p$ is symmetric in the arguments $\lambda_1,\ldots,\lambda_k$ for each $k$, 
$p$ is normalised by the condition $p(1)=1$, and $p$ satisfies 
the addition rule
\begin{equation}\label{add-rule}
p(\lambda)=\sum_{\mu: \,\mu\searrow\lambda}p(\mu)
\end{equation}
where the sum is over compositions $\mu$ derived from $\lambda$ by either increasing
a part by one or by appending $1$ at the end of the sequence $\lambda$.
For instance, 
if $\lambda=(3,2,2)$, $\mu$ assumes the values $(4,2,2)$, $(3,3,2)$, $(3,2,3)$ and $(3,2,2,1)$,
and (\ref{add-rule}) specialises to
$$p(3,2,2)=p(4,2,2)+2 p(3,3,2)+p(3,2,2,1).$$ 
A function $p$ with these properties is known as 
an {\it exchangeable partition probability function} (EPPF).
Such a function uniquely determines the probability law of a corresponding exchangeable 
random partition $\Pi$.
\par According to a Kingman's {\em paintbox representation} \cite{jp.bmpart, gnedin97, gnedinp03}, 
every such exchangeable partition has the same distribution as $\Pi$ constructed from
some random closed set $Z\subset [0,1]$, as follows:
let $u_1, u_2, \ldots$ be independent uniform $[0,1]$ variables, independent of $Z$,
and let distinct integers $i$ and $j$ belong to the same block of $\Pi$ if and only 
$u_i$ and $u_j$ fall in the same open interval component of $[0,1] \backslash Z$.

\par A distinguished class of exchangeable partitions is the two-parameter family, with EPPF
\begin{equation}\label{2par}
p_{\alpha,\theta}(\lambda_1,\ldots,\lambda_k):={(\theta+\alpha)_{k-1\uparrow \alpha}\over 
(\theta+1)_{n-1\uparrow }}\prod_{j=1}^k (1-\alpha)_{\lambda_j-1\uparrow }
\end{equation}
where $n=\Sigma\, \lambda_j$ and 
$$(x)_{m\uparrow\beta}:=\prod_{j=1}^m (x+(j-1)\beta)\,,~~~~~~(x)_{m\uparrow}=(x)_{m\uparrow 1}$$
are rising factorials,
with the convention that $(x)_{0 \uparrow \beta } := 1$.
Possible values of the parameters $(\alpha,\theta)$ are 
either $-\infty\leq \alpha<0$ and $\theta=m|\alpha|$ for some $m=1,2,\ldots,\infty\,$;
or $0\leq \alpha\leq 1$ and $\theta\geq -\alpha$,
with a proper understanding of (\ref{2par}) in some limiting cases.
See \cite{csp} for detailed exposition of the general theory of exchangeable partitions 
and features of the $(\alpha,\theta)$ family.

\par In this paper we are interested in a special class of {\it Gibbs partitions}, 
which generalise (\ref{2par}) as follows:
\vskip0.5cm
\begin{definition}{\rm An 
exchangeable random partition $\Pi$ of the set of natural numbers 
is said to be of {\it Gibbs form} if for some nonnegative weights $\W=(\W_j)$ and $\V=(\V_{n,k})$
the EPPF of $\Pi$ satisfies 
\begin{equation}\label{gibbs}
p(\lambda_1,\ldots,\lambda_k)=\, \V_{n,k}\prod_{j=1}^k \W_{\lambda_j}\,
\end{equation}
for all $1\leq k\leq n$ and all compositions $(\lambda_1,\ldots,\lambda_k)$ of $n$.
}
\end{definition}
\noindent

\par For fixed $n$ we can choose arbitrary nonnegative weights $\V_1,\ldots,\V_n$ and $\W_1,\ldots,\W_n$ which are
not identically zero,
and use (\ref{gibbs}) to define a random partition of $[n]$ by setting
$\V_{n,k}=\V_k/c_n$ for $c_n$ a suitable normalisation constant 
(see \cite{Ve} for another version of the Gibbs formalism).
The block sizes of such a Gibbs partition can be realised by Kolchin's model, that is 
identified with the collection of terms of a random sum $S=X_1+\ldots+X_K$ conditioned on 
$S=n$, with independent identically distributed $X_1, X_2,\ldots$, independent of $K$.
For integer weights this is the distribution on partitions of $[n]$ induced by components of 
a random composite structure built over partitions of $[n]$, 
when there are $\W_j$ possible configurations associated with every subset of $[n]$ with $j$ elements, $\V_k$ possible 
configurations 
associated with every collection of $k$ blocks, and a uniform distribution is assigned to all possible composite structures
subject to these constraints.
For instance, if 
$\W_j=(j-1)!$ and $\V_k=\theta^k$ (with $\theta\in {\mathbb N}$), the product $\V_k \W_{\lambda_1}\cdots \W_{\lambda_k}$
counts the number of coloured permutations of
$[n]$ with cycle sizes $(\lambda_1,\ldots,\lambda_k)$ and one of $\theta$ possible colours
assigned to each of the cycles: then (\ref{gibbs}) reduces to (\ref{2par}) with $\alpha=0,$ 
$\theta>0$.
So in this case, there is an infinite exchangeble partition $\Pi$ whose restrictions $\Pi_n$
are all of the Gibbs form \re{gibbs}.
Many other combinatorially interesting examples of Gibbs partitions $\Pi_n$ can be given, 
using the prescription \re{gibbs} for each fixed $n$: see for instance \cite{abt, csp}. But 
typically the distributions of these combinatorially defined $\Pi_n$  are not consistent as 
$n$ varies, so they are not realisable as the sequence of restrictions to $[n]$ of
an infinite Gibbs partition.

\par The special case of $\V$-weights representable as ratios
$\V_{k,n}=\V_k/c_n$ 
was studied 
by Kerov \cite{kerov.cra} in the framework of Kolchin's model.
In this case one assumes a {\it single} infinite sequence of weights $(\V_k)$ and the $c_n$'s appear as normalisation 
constants. Kerov \cite{kerov.cra} established that the Gibbs partitions of this type 
are precisely the members of the two-parameter family (\ref{2par}).

\par We will show that in the more general setting (\ref{gibbs}), allowing an arbitrary
triangular array $\V_{k,n}$, the $\W$-weights must be still as in 
(\ref{2par}), with a single parameter $\alpha\in [-\infty, 1]$ defining the {\it type} of Gibbs 
partition. 
The case $\alpha=1$ is trivial. 
For each nontrivial type
$\alpha<1$ the set of all possible $\V$-weights is an infinite simplex ${\cal V}_\alpha$.
We identify the extreme elements of ${\cal V}_\alpha$ by solving a
boundary problem for an instance of the generalised Stirling triangle, as introduced in 
another paper by Kerov \cite{kerov.af} (also see \cite[Chapter I]{kerov.diss}).
It turns that the nature of the extremal set depends substantially on the type.
According to our main result, stated more formally in Theorem \ref{mainth}, there are three
qualitatively different ranges of $\alpha$.
For $\alpha \in [-\infty,0[$ the extremal set is discrete and corresponds to the 
members of the $(\alpha,\theta)$-family. 
For $\alpha=0$ this set is continuous and still corresponds to the members of the $(\alpha,\theta)$-family
(Ewens' partitions). 
For $\alpha \in ]0,1[$ the $(\alpha,\theta)$-partitions are not extreme, rather the 
extremes of ${\cal V}_\alpha$ comprise 
a continuous $(\alpha|s)$-family (with parameter $s\in [0,\infty]$) which appears by conditioning the $(\alpha,\theta)$-partitions on the asymptotics 
of the number of blocks. In \cite{jp.pk} 
the $(\alpha|s)$-partitions were derived from their Kingman's representation, with the random
closed set $Z$ being the scaled range of an $\alpha$-stable subordinator 
conditioned on its value at a fixed time. This identification of extreme elements of ${\cal V}_\alpha$ for $0\leq\alpha<1$
was indicated without proof in \cite[Theorem 8]{jp.pk}.

\section{Some basic results}

To define an exchangeable partition the weights in (\ref{gibbs}) are normalised by the condition
$$ 
\sum_{k=1}^n \V_{n,k} B_{n,k}(\W)=1~~~{\rm for~~}n=1,2,\ldots
$$
where $B_{n,k}$ is a partial Bell polynomial in the variables $\W=(\W_1, \W_2,\ldots)$, 
\begin{equation}\label{bell}
B_{n,k}(\W):=\sum_{\{A_1,\ldots,A_k\}}\prod_{j=1}^k \W_{|A_j|}={n!\over k!} 
\sum_{(\lambda_1,\ldots,\lambda_k)}\prod_{i=1}^k {\W_{\lambda_j}\over \lambda_j!}\,,
\end{equation}
where the first sum expands over all partitions of $[n]$ into $k$ blocks, and the second 
over all compositions of $n$ with $k$ parts.
Observe that there is a redundancy in the possible values of parameters:
the distribution of the Gibbs partition is unaffected by simultaneous substitutions when either
$\W_j\to \gamma^j\W_j$ and $\V_{n,k}\to \gamma^{-n}\V_{n,k}$, or 
$\W_j\to \gamma \W_j$ and $\V_{n,k}\to \gamma^{-k}\V_{n,k}$
for $\gamma>0$. 
Throughout we assume the 
normalisation $\V_{1,1}=\W_1=1$. Granted the normalisation and excluding the trivial case
of a partition $\Pi$ with only singleton blocks
the ambiguity amounts to the geometric tilting
$\W_j\to \gamma^{j-1}\W_j$ and $\V_{n,k}\to\gamma^{k-n}\V_{n,k}$.
\par Our starting point is the following elementary lemma. 
\begin{lemma}\label{basic}
The weights $(\W_j)$ and $(\V_{n,k})$ with $\W_1=\V_{1,1}$ 
define a partition of Gibbs form if and only if for some $b\geq 0$ and $a\leq b\,$ the following two 
conditions are satisfied:
\begin{itemize}
\item[{\rm (i)}] 
\begin{equation}\label{wab}
\W_j=(b-a)_{j-1\uparrow b}\,,~~j=1,2,\ldots
\end{equation}
\item[{\rm (ii)}] the $\V_{n,k}$ satisfy the recursion
\begin{equation}\label{rec}
\V_{n,k}=(bn-ak)\V_{n+1,k}+\V_{n+1,k+1}\,,~~~~~~1\leq k\leq n\,.
\end{equation}
\end{itemize} 
\end{lemma}
\noindent
{\it Proof.} The trivial
singleton partition of ${\mathbb N}$
is of the Gibbs form with either
$\W_j=1(j=1)$ or $\V_{n,k}=1(k=n)$. 
Excluding the singleton partition,
the Gibbs prescription forces
$\W_j>0$ for all $j=1,2,\ldots$, because $p(n)=\V_{n,1} \W_n$ and 
$p(n)>0$ for all $n$ (as follows from Kingman's representation).

The only Gibbs partition with $\V_{n,2}=0$ for some $n\geq 2$ is the trivial one-block partition, 
in which case the conclusion is obvious. 
Excluding also this trivial case we may assume that
$\V_{n,2}>0$ for all $n$ and $\W_j>0$ for all $j$.
Introducing $r_j=\W_j/\W_{j+1}$ we find then
that (\ref{add-rule}) amounts to
$$\V_{n,k}=\V_{n+1,k}\sum_{j=1}^k r_{\lambda_j}+\V_{n+1,k+1}$$
for all compositions of $n$ with $k$ parts. Applying this for $k=2$ and using $\V_{n,2}>0$ we see that $r_i+r_j$ depends only 
on $i+j$, hence $r_{j+1}-r_j$ is constant and therefore $(r_j)$ is an arithmetic sequence
$$
r_j=bj-a\,,~~~j=1,2,\ldots
$$
where necessarily $b\geq 0$ and $a<b$ to ensure $r_j>0$ (as entailed by $\W_j>0$).
Now, from $\W_j=r_1\cdots r_{j-1}$ we obtain (\ref{wab}), and (\ref{rec}) follows because
$$\sum_{j=1}^k r_{\lambda_j}=b\sum_{j=1}^k \lambda_j-ak\,.$$
Inverting the argument we see that $p$ defined by (\ref{wab}), (\ref{rec}) and (\ref{gibbs})
satisfies the addition rule (\ref{add-rule}).
\endpf
\noindent

\par Two parameters $a$ and $b$ may be reduced by the geometric tilting 
to a single parameter $\alpha\in \,[-\infty,1]$,
corresponding to the sequence of $\W$-weights 
\begin{equation}\label{weights}
\W_j=\left\{\begin{array}{l}
(1-\alpha)_{j-1\uparrow}~~~{\rm for~~} -\infty<\alpha\leq 1\,,\\
~~~~1~~~~~~~~~~~~\,{\rm for~~~}
\alpha=-\infty\,.
\end{array}\right.\,
\end{equation}
The case $\alpha=1$ corresponds to the trivial singleton partition and will be excluded from
further consideration.

\begin{definition}{\rm 
For $\alpha<1$ 
let ${\cal P}_\alpha$ be the set of all distributions of infinite partitions $\Pi$
of {\it type} $\alpha$, 
whose EPPF
$p$ is of the Gibbs form (\ref{gibbs}) with these $\W$-weights (\ref{weights}), 
and
let ${\cal V}_\alpha$ be the set of nonnegative solutions $\V=(\V_{n,k})$ to the backward recursion
\begin{equation}\label{rec1}
\V_{n,k}=
\gamma_{n,k}
\V_{n+1,k}+\V_{n+1,k+1}\,,
\end{equation}
with $\V_{1,1}=1$, where the coefficient for $1\leq k\leq n$ is given by 

\begin{equation}\label{gamma}
\gamma_{n,k}=\left\{\begin{array}{l}
n-\alpha k~~~{\rm for~~} -\infty<\alpha<1\,,\\
~~~~k~~~~~~\,{\rm for~~~}
\alpha=-\infty\, .
\end{array}\right.\,
\end{equation}
}
\end{definition}

Lemma \ref{basic} establishes an affine 
bijection between these two convex sets ${\cal P}_\alpha$ and ${\cal V}_\alpha$, 
hence also a bijection between the sets of their extreme points. 
To spell this out, each probability distribution $\prob$ of $\Pi$ with 
$\prob \in {\cal P}_\alpha$ induces a distribution of $K_n$, the number of blocks of
$\Pi_n$, according to the formula
$$
\prob(K_n=k)=\V_{n,k}B_{n,k}(\W)\,,
$$
obtained by summation of \re{gibbs}
over all partitions of $[n]$ with $k$ blocks. On the other hand,
the conditional distribution of $\Pi_n$ given the number of blocks is
\begin{equation}\label{cond}
{\mathbb P}(\Pi_n=\{A_1,\ldots,A_k\}|K_n=k)={\prod_{j=1}^k \W_{|A_j|}\over B_{n,k}(\W)} .
\end{equation} 
Thus the distribution $\prob$ of $\Pi$ determines the weights $\V_{n,k}$, and vice versa.
Moreover, the weak topology on ${\cal P}_\alpha$, defined by pointwise convergence of
EPPF's, corresponds in ${\cal V}_\alpha$ to convergence of 
the $\V_{n,k}$ for all $1 \le k \le n$.

\par Observe that (\ref{cond}) does not involve the $\V$-weights. Thus 
for Gibbs partitions of a given type $\alpha$,
the sequence of block counts $(K_n)$ is a sequence of sufficient statistics for $(\Pi_n)$.
In particular, for $(\alpha,\theta)$-partitions the sequence $(K_n)$ is a Markov chain whose
time-reversed transition probabilities are the same for all $\theta$. 

\par By general theory of sufficient statistics and extreme points \cite{dy78, diaconis},
each $p\in {\cal P}_\alpha$ can be uniquely represented as a convex mixture of the extreme elements of ${\cal P}_\alpha$, and the same can be said of ${\cal V}_\alpha$.

\par Let $\phi(\alpha,\theta):=(\phi_{n,k}(\alpha,\theta), 1\leq k\leq n)$ denote the particular
sequence of $\V$-weights appearing in the two-parameter formula (\ref{2par}),
that is 
\eq
\label{phink}
\phi_{n,k}(\alpha,\theta)={({\theta+\alpha})_{k-1\uparrow\alpha}\over (\theta+1)_{n-1\uparrow}}\,\qquad {\rm for~~}
\alpha\neq - \infty
\en
where $\theta\geq \alpha$ for $0\leq\alpha<1$, and $\theta=-|\alpha|m$, with $m=1,2,\ldots,\infty$ for $\alpha<0$. 
In the case $\alpha=-\infty$ the formula 
\eq\label{phinf}
\phi_{n,k}(-\infty,m\infty)={(m)_{k\downarrow}\over m^n}\qquad
\en
is the $\alpha\to-\infty$ limit 
for $m=1,2,\ldots,\infty$, with $m=\infty$ corresponding to the trivial singleton partition.
We have $\phi(\theta,\alpha)\in {\cal V}_\alpha$, as can be readily checked by algebra.
A characteristic property of this class of solutions is summarised in the following corollary: 

\begin{corollary}{\rm \cite[Theorem 7.1]{kerov.cra}} Members of the two-parameter family, with EPPF given by
{\rm (\ref{2par})}, are the only partitions of the Gibbs form with $\V$-weights representable as ratios
$\V_{n,k}=\V_k/c_n$.
\end{corollary}
{\it Proof.} Assuming $\V_{\varkappa}>0$ for some 
$\varkappa\geq 2$ exchangeability implies $\V_{k}>0$ for $k\leq\varkappa$ and the recursion
(\ref{rec}) becomes 
$$
{\V_{k+1}\over \V_k}-ak={c_{n+1}\over c_n}-bn\,,~~~1\leq k\leq \varkappa\,.
$$
Since the left side does not involve $n$ and the right side does not involve $k$, their common value 
is a constant, say $t$. Then $\V_k=(t+a)_{k-1\uparrow a}$, and because $\V_2>0$ we have $t>-a$. In the case
$a<0$ the requirement $\V_k\geq 0$ forces $t=-ma$ for some integer $m\geq 1$. By a similar argument
$c_n=(t+b)_{n-1\uparrow b}$. The form (\ref{2par}) follows by redundancy. \endpf

\section{The boundary problem}

To embed our discussion of extremes
in a wider context let $\cal V$ be the convex set of nonnegative solutions $\V=(\V_{n,k})$ to the
recursion 
\begin{equation}\label{pascal}
\V_{1,1} = 1;~~~~
\V_{n,k}=\gamma_{n,k}\V_{n+1,k}+\delta_{n,k}\V_{n+1,k+1}\,, ~~~~~(1\leq k\leq n)
\end{equation}
where the coefficients $\gamma_{n,k}$ and $\delta_{n,k}$ form two arbitrary triangular
arrays of non-negative numbers.
This recursion 
is associated with 
a generalised Pascal triangle, that is 
an infinite directed graph $G$ with vertex set 
$\{(n,k):1\leq k\leq n\}$, such that each 
vertex $(n,k)$ has two immediate successors $(n+1,k)$ and $(n+1,k+1)$, the
multiplicities of the outgoing edges being $\gamma_{n,k}$ and $\delta_{n,k}$, respectively.
For a directed path connecting the root $(1,1)$ and a node $(n,k)$ we define the {\it weight} of the path 
to be the product of the multiplicities
along the path, 
and let the {\it dimension} $d^{n,k}$ be the sum of weights
of all such paths, with the convention 
$d^{1,1}=1$. 
The dimension is a unique solution to the forward recursion
\begin{equation}\label{dpascal}
d^{n+1,k}=\delta_{n,k-1}d^{n,k-1}+\gamma_{n,k}d^{n,k}
\end{equation}
where $\delta_{n,0}=0$. 
In the case of the standard Pascal triangle with $\gamma_{n,k}\equiv 1, \,\delta_{n,k}\equiv 1$
this reduces by an obvious shift of indices 
to the familiar recursion for the binomial coefficients.

\par Note that each path in $G$ from $(1,1)$ to $(n,k)$ may be written as a sequence $k_1, k_2,\ldots,k_{n}$
with $k_1=1$, $k_n=k$ and $k_{j+1}-k_j\in \{0,1\}$. 
Consider now a random process $(K_n)$ with $1\leq K_n\leq n$ 
such that for every path from $(1,1)$ to $(n,k)$ 
the conditional probability 
$${\mathbb P}(K_1=1,K_2=k_2,\ldots,K_{n-1}=k_{n-1}|K_n=k)$$
equals the weight of this path divided by $d^{n,k}$. 
That is to say, the process $(K_n)$ is Markovian with co-transition probabilities
$$
{\mathbb P}(K_{n-1} = j | K_n=k) = {d_{n-1,j}\, q_{n-1}(j,k) \over d_{n,k} }
$$
where
\begin{equation}\label{cotrans}
q_{n-1}(j,k) = \left\{\begin{array}{l}
\gamma_{n-1,j} {\rm ~~if ~~} j = k \\
\delta_{n-1,j} {\rm ~~if ~~} j = k - 1\\
~0 {\rm~~ ~~~~else }\,.
\end{array}\right.\,
\end{equation}
The formula
$${\mathbb P}_{\V} (K_n=k)={\V}_{n,k}d_{n,k}$$ 
establishes a bijection between $\cal V$ 
and the set of laws of such Markov chains. So we identify
$\cal V$ with this set of Markovian laws.

\par By some well known general theory 
\cite{aldous,diaconis,dy78,kemeny,kerov.young}), each extreme law in ${\cal V}$
can be represented as a weak limit of the
conditional laws for $(K_n)$
given $K_\nu=\varkappa_\nu$ for $\nu\to\infty$ and some sequence $(\varkappa_\nu, \nu=1,2,\ldots)$. 
Following \cite{koo}
we will call the set of these limit laws the {\it boundary} of $G$ 
(sometimes also called `Martin boundary' or `maximal boundary' or `the set of Boltzmann laws').

\par More explicitly, 
extending the above definitions,
let the weight of a path in $G$ connecting $(n,k)$ and $(\nu,\varkappa)$
be the product of multiplicities along this path, and define
the extended dimension
$d_{n,k}^{\nu,\varkappa}$ to be the total weight of paths in $G$ connecting
$(n,k)$ and $(\nu,\varkappa)$ (which is zero unless $ k\leq \varkappa \leq\nu\,$ and $\nu\geq n$).
Thus 
$d_{1,1}^{n,k}=d^{n,k}$ and $d_{n,k}^{n,k}=1$. 
For each chain directed by some $\V \in {\cal V}$ the conditional law of $(K_1,\ldots,K_{n-1})$ given
$K_\nu=\varkappa$ (for $\nu>n$) is the same, and is determined by 

\eq\label{Vnk}
{\mathbb P}(K_n=k\,|\,K_\nu=\varkappa)=\V^{\nu,\varkappa}_{n,k}\, d^{n,k}\,,~~~
{\rm where~~} V^{\nu,\varkappa}_{n,k}:=\,{d_{n,k}^{\nu,\varkappa}\over d^{\nu,\varkappa}}\,.
\en
Clearly, for $(\nu,\varkappa)$ fixed, $\V^{\nu,\varkappa}_{n,k}$ satisfies 
(\ref{pascal}) for $n<\nu$, thus if $\V^{\nu,\varkappa_\nu}_{\bullet,\bullet}\,$ converge as $\nu\to\infty$ 
along some infinite path
$(\varkappa_\nu)$ then the limit is certainly in ${\cal V}$. 
The infinite paths which induce the limits are called {\it regular}, and 
the set of such limit elements of $\cal V$ is the boundary of $G$.
For regular path $(\varkappa_{\nu})$ we say that the boundary element $\V=\lim_{\nu\to\infty}\V^{\nu,\varkappa_\nu}$ 
and the corresponding law ${\mathbb P}_{\V}$ are {\it induced} by the path. 

\par The next lemma is adapted from \cite{aldous, diaconis, dy78, kemeny}.

\begin{lemma}\label{gen} Identifying the elements of ${\cal V}$ with the laws ${\mathbb P}_{\V}$ for Markov chain 
$(K_n)$ we have: 
\begin{itemize}
\item[\rm (i)] each extreme element $\V\in {\cal V}$ belongs to the boundary of $G$, that is
may be represented as a limit of 
the functions $\V^{\nu,\varkappa_\nu}_{\bullet,\bullet}\,$ along some regular path $(\varkappa_\nu)$,
\item[\rm (ii)] for every $\V\in {\cal V}$, 
under 
${\mathbb P}_{\V}$ almost all paths of $(K_n)$ are regular, 
\item[\rm (iii)] a solution $\V\in {\cal V}$ is extreme iff the set of regular paths which induce 
$V$ has $\prob_{\V}$-probability one. 
\end{itemize}
\end{lemma}
\vskip0.5cm
\noindent
{\bf Example.}
The instance of the boundary problem for the standard Pascal triangle
has been treated by many authors.
In this case it is more convenient 
to label the nodes by nonnegative integers $\{(n,k), 0\leq k\leq n\}$.
The corresponding chains $0 \le K_n \le n$ are those whose increments
$K_{n+1} - K_n$ are exchangeable random variables with values in $\{0,1\}$.
The dimension function is given by the binomial 
coefficients 
$d^{n,k}={n\choose k}\,,~~d_{n,k}^{\nu,\varkappa}={\nu-n\choose\varkappa- k}$. 
A path $(\varkappa_\nu)$ is regular if and
only if there is a limit $\varkappa_\nu/\nu\to s$ for some $s\in [0,1]$, which corresponds to a
boundary element $\V(s)$ with $\V_{n,k}(s)=s^{k}(1-s)^{n-k}$, hence the boundary is homeomorphic to $[0,1]$. 
This is de Finetti's representation of infinite exchangeable sequences of zeros and ones.
Since 
$K_n$ under ${\mathbb P}_{\V(s)}$ is the number of successes
in a series of $n$ Bernoulli trials with success probability $s$, the law of large numbers ensures
${\mathbb P}_{\V(s)}(K_n/n\to s)=1$. Hence each $\V(s)$ is extreme 
by Lemma \ref{gen} (iii).
In fact, to ensure regularity of a path we only need to check the convergence of
$\V_{n,0}^{\nu,\varkappa_\nu}$ for each $n$, because the bivariate array $(\V_{n,k})$ satisfies 
the backward Pascal recursion if each $\V_{n,k}$ is a finite difference of the sequence
$(\V_{n,0})$. It follows that a sequence $(\V_{n,0})$ with $\V_{n,0}=1$ is representable 
as a convex mixture of functions $\V_{n,0}(s)=(1-s)^n, \,s\in [0,1]$ if and only if 
the associated array $(\V_{n,k})$ is nonnegative, in which case such representation is unique.
The last assertion is widely known as the resolution of
the Hausdorff problem of moments.
\vskip0.5cm
\noindent
In general, however, 
the set of extremes 
(sometimes called the `minimal' boundary) may be smaller than the boundary.
This kind of pathology is illustrated by the following example. 
\vskip0.5cm
\noindent
{\bf Example} 
Consider a graph $G$ with the following sets of nodes and edges.
Level $0$ has a single node $\emptyset$, which is the root of $G$.
Level $1$ has two nodes $a_1$ and $c_1$, 
and the root $\emptyset$ is 
connected to the nodes $a_1$ and $c_1$.
Level $2$ has three nodes $a_2, b_2, c_2$, so that $a_1$ is connected to $a_2$ and $b_2$, while $c_1$ is connected
to $b_2$ and $c_2$. On each further level $n>2$ there are exactly $3$ nodes $a_n, b_n, c_n$. Node
$a_{n-1}$ is connected to $a_n$ and $b_n$, node $b_{n-1}$ is connected only to $b_n$, and $c_{n-1}$ is connected to $b_n$ and $c_n$.
There are no other edges and the
edges just described all have multiplicity $1$.
\par Every infinite path in $G$ starting at $\emptyset$ is regular.
The boundary of $G$ consists of three elements $V_a, V_b$ and $V_c$, induced by the paths 
$a=(\emptyset, a_1,a_2,\ldots)$, $b=(\emptyset, b_1,b_2,\ldots)$, 
and $c=(\emptyset, c_1,c_2,\ldots)$, respectively.
Clearly, $V_a$ is a unit mass at $a$ and
$V_c$ is a unit mass at $c$.
Observe that for $\nu>n\geq 2$ there are $\nu-n$ paths from $a_n$ to $b_\nu$, the same number
$\nu-n$ of paths from $c_n$ to $b_\nu$, and there is only one path connecting $b_n$ and $b_\nu$.
Sending $\nu\to\infty$ we see that 
$V_b$ is the mixture $V_b=V_a/2+ V_c/2$.

\par It follows that the boundary $\{V_a, V_b, V_c\}$ is larger than the set of extremes $\{V_a, V_c\}$.
The distribution $V_b$ violates the condition in Lemma \ref{gen} (iii): though
$V_b$ can be induced by many paths (unlike $V_a$ and $V_c$), the set of these paths has $V_b$-probability zero.

\vskip0.5cm
\par Further examples may be related to other classical number triangles and their generalisations 
(as in 
\cite{labelle, regev}) although explicit results on the boundary problem are scarce.
For later application we record some useful general tools.

\vskip0.5cm

\par Obviously, $\V_{n,k}=0$ implies $\V_{\nu,\varkappa}=0$ for $\nu\geq n,$ $\varkappa\geq k$.
In particular, the trivial solution with ${\mathbb P}_{\V} (K_n= n)=1$ is characterised by $\V_{2,1}=0$,
and another trivial solution with ${\mathbb P}_{\V}(K_n=1)$ is characterised by $\V_{2,2}=0$.
Both trivial solutions are extreme.
\begin{lemma}\label{reg}
A path $(\varkappa_\nu)$ is regular if and only if 
$\V_{n,1}^{\nu,\varkappa_\nu}$ converge as $\nu\to\infty$ for each $n$. 
Every path with $\V_{2,1}^{\nu,\varkappa_\nu}\to 0$ is regular and induces
the trivial solution with $\V_{2,1}=0$.
\end{lemma}
{\it Proof.} Because
$\V_{n+1,k+1}=(\V_{n,k}-\gamma_{n,k}\V_{n+1,k})/\delta_{n,k}\,$,
a double
induction, first in $k$ and then in $n$, shows that $v\in {\cal V}$ is uniquely determined by the 
entries $(\V_{n,1})$. 
\endpf

The next lemma expresses a well known {\em stochastic monotonicity} 
property of the kind of inhomogenous positive integer-valued
Markov chains involved here.
We indicate an algebraic proof, but it can also be derived probabilistically 
by a coupling argument. See \cite{aldous83} and papers cited there.

\begin{lemma}\label{monotone} For $\nu\geq n$ fixed,
$V_{n,1}^{\nu,\varkappa}$ is nonincreasing in $\varkappa$.
\end{lemma}
{\it Proof.} The proof is by induction in $\nu$. 
Suppose the claim is true for some $\nu$, then for fixed $1\leq\varkappa\leq n$ 
and nonnegative $\alpha, \beta, \gamma, \delta$ with $\alpha+\beta=1,\,\gamma+\delta=1\,$
we have
$$\alpha V_{n,1}^{\nu,\varkappa-1} +\beta V_{n,1}^{\nu,\varkappa}\geq
V_{n,1}^{\nu,\varkappa}\geq \gamma V_{n,1}^{\nu,\varkappa}+\delta V_{n,1}^{\nu,\varkappa+1}$$
(where $ V_{n,1}^{\nu,0}= V_{n,1}^{\nu,\nu+1}=0$).
For a suitable choice of $\alpha,\beta,\gamma,\delta$ the left side 
of the inequality equals $d^{n,1} V_{n,1}^{\nu+1,\varkappa}$ while the right side equals
$d^{n,1} V_{n,1}^{\nu+1,\varkappa+1}$, as follows readily from 
(\ref{Vnk}). The induction step follows. 
\endpf

\begin{lemma}\label{disc}
Suppose for $m=1,2,\ldots$ there are solutions $\V(m)\in {\cal V}$ such that 
$\V_{n,m} (m)\,d^{n,m}\to 1$ as $n\to\infty$, 
then each $\V(m)$ is extreme and satisfies
\begin{equation}\label{pvm}
{\mathbb P}_{\V(m)}\left(\lim_{n\to\infty} K_n= m\right)=1\,.
\end{equation}
If above that $\V_{2,1}(m)\to 0$ as $m\to\infty$
then $\V(m)$ converges to the trivial law $\V(\infty)$ with 
${\mathbb P}_{\V(\infty)}(K_n=n)=1$,
and in this case the set of extreme elements of $\cal V$ is
$\{\V(1), \V(2),\ldots, \V(\infty)\}$. 
\end{lemma}
{\it Proof.} All paths $(\varkappa_\nu)$ are nondecreasing, thus ${\mathbb P}_{\V(m)}(K_n>m)=0$, and because
$${\mathbb P}_{\V(m)}\left(\lim_{n\to\infty}K_n= m\right)\geq {\mathbb P}_{\V(m)}(K_n=m)=\V_{n,m} (m)\,d^{n,m}\to 1
\,,~~~{\rm as~} n\to\infty$$
we have (\ref{pvm}). 
Easily from Lemma \ref{gen}, $V(m)$ is extreme and can be induced by arbitrary 
path with $\varkappa_\nu=m$ for large enough $\nu$.

\par Now let $\varkappa_\nu\uparrow\infty$.
By Lemma \ref{monotone} and the above argument 
we have for $\nu\to\infty$
$$\V_{2,1}^{\nu,\varkappa_\nu}<\V_{2,1}^{\nu,m}\to \V_{2,1}(m)\,,$$
hence letting $m\to\infty$ and invoking Lemma \ref{reg} shows that $(\varkappa_\nu)$ induces
$\V(\infty)$.
Since every path has either finite or infinite limit, every path is regular and the
list of extremes is complete.
\endpf
\noindent


\vskip0.5cm
\noindent
{\bf Example} 
A discrete family of solutions with the properties as in Lemma \ref{disc}
exists for a graph called 
`the $q$-Pascal triangle'. 
The set of nodes of the graph is $\{(n,j),0\leq j\leq n\}$,
the multiplicities are $\gamma_{n,k}=q^{k-1},\,\delta_{n,k}=1$, and
the dimension function is given by the $q$-binomial coefficients.
The boundary has been determined in 
\cite{kerov.af, olshanski2001}.

\vskip0.5cm
\par Our main tool for identifying the boundary in the continuous
case is the following lemma. Compare with \cite{jp.def} where 
the same method is applied to obtain a different generalisation of
de Finetti's theorem for sequences of zeros and ones, and see
\cite{aldous83} for another closely related setting.

\vskip0.5cm

\begin{lemma}\label{cont} 
Suppose there is a sequence of positive constants $(c_n)$
with $c_n\to\infty$,
and for each $s\in [0,\infty]$ there is
a solution $V(s)\in {\cal V}$ which satisfies 
\eq\label{lln}
{\mathbb P}_{\V(s)}\left(\lim_{\nu\to\infty}K_\nu/c_\nu= s\right)=1\,.
\en
Suppose 
the mapping $s\mapsto \V(s)$ is a continuous injection from $[0,\infty]$ to ${\cal V}$ with
$0$ and $\infty$ corresponding to the trivial solutions
$${\mathbb P}_{\V(0)}(K_n=1)=1\,, ~~~{\mathbb P}_{\V(\infty)}(K_n=n)=1\,.$$
Then 
\begin{enumerate}
\item[\rm (i)] a path $(\varkappa_\nu)$ is regular if and only if $\,\lim_{\nu\to\infty}\varkappa_\nu/c_\nu=s$ for some 
$s\in [0,\infty]$, in which case $(\varkappa_\nu)$ induces $\V(s)$,
\item[\rm (ii)] $\{\V(s), s\in [0,\infty]\}$ is the set of extreme elements 
of ${\cal V}$. 
\end{enumerate}
\end{lemma}

\noindent
{\it Proof.} (i) 
Let $(\varkappa_\nu)$ be a path with
$\varkappa_\nu/c_\nu\to s$ for some $0<s<\infty$.
Using the fact that for all $V\in {\cal V}$ the co-transition probabilities are the same,
and exploiting the monotonicity, as in Lemma \ref{monotone}, we can 
squeeze 
$${\mathbb P}_{V(s-2\epsilon)}(K_n=1|K_\nu/c_\nu<s-\epsilon)>d^{n,1} V_{n,1}^{\nu,\varkappa_\nu}>
{\mathbb P}_{V(s+2\epsilon)}(K_n=1|K_\nu/c_\nu >s+\epsilon)$$
for $\epsilon<s/2$ and $\nu$ sufficiently large.
From this and the assumption (\ref{lln}) we derive
$$V_{n,1}(s-2\epsilon)\geq V_{n,1}^{\nu,\varkappa_\nu}\geq V(s+2\epsilon)$$
for large $\nu$. Letting $\epsilon\to 0$ and using the assumed continuity we conclude that $(\varkappa_\nu)$ is regular 
and induces $V(s)$. The cases $s=0$ and $s=\infty$ are treated similarly.
By the same argument, a path $(\varkappa_\nu)$ 
cannot be regular if $\varkappa_\nu/c_\nu$ has distinct subsequence limits. 
\par (ii) Follows from (i), (\ref{lln}) and Lemma \ref{gen} (iii).
\endpf

\noindent
The lemma is designed to cover normalisations $c_n=o(n)$. Compare this with the 
standard Pascal triangle, where we assume the scaling by $c_n=n$ that leads to parameterisation 
of the boundary by $[0,1]$.

\section{Stirling triangles}

A {\it generalised Stirling triangle}, as introduced by 
Kerov \cite{kerov.af}, is a generalised Pascal graph $G$ with 
multiplicities of the form 
$\gamma_{n,k}=b_n+a_k\,,\,$ and $\delta_{n,k}=1$.
We will consider the boundary problem in the special case with coefficients (\ref{gamma})
where $\alpha\in [-\infty,1[\,.$
In this case the dimension 
$d^{n,k}$ is a generalised Stirling number 
$\left[\!\!\begin{array}{c} n\\k \end{array}\!\!\right]_\alpha$ which may be defined 
in many ways. For instance, it is determined
by the recursion
(\ref{dpascal}), or by specialising 
the Bell polynomial $B_{n,k}$ for weights (\ref{weights}),
or as 
the connection coefficient in 
$$(x)_{\uparrow}=\sum_{k=1}^n \left[\!\! \begin{array}{c} n\\k \end{array}\!\!\right]_\alpha (x)_{k\uparrow\alpha}\,
{\rm ~~~(for}~~\alpha\neq-\infty)\,,
$$
or as coefficient at $x^{n}$ in the series expansion of
\begin{equation}\label{genfun}
{n!\over \alpha^k\,k!}\,(1-(1-x)^\alpha)^k\,,\qquad (\alpha\neq \,0,-\infty)\,.
\end{equation}

\noindent
For $\alpha=-\infty$ these are the Stirling numbers of the second kind, for $\alpha=0$
the signless Stirling numbers of the first kind and for $\alpha=-1$ the Lah numbers. 
Alternatively, by the definition of dimension $d^{n,k}$ as the sum of weights we obtain 
\begin{eqnarray}\label{St}
\left[\!\! \begin{array}{c} n\\ k \end{array}\!\!\right]_\alpha &=&
\sum_{1=n_0< n_1<\ldots<n_{k-1}<n_k=n}~\prod_{j=1}^k~\prod_{n_{j-1}<\nu<n_j}(\nu-\alpha j)\,,
~~~{\rm for~~}\alpha > - \infty
\\
\left[\!\! \begin{array}{c} n\\ k \end{array}\!\!\right]_{-\infty}&=&
\sum_{1=n_0< n_1<\ldots<n_{k-1}<n_k=n}~\prod_{j=1}^k~j^{n_j-n_{j-1}-1}\,\,.
~~~~~~~~~~~~~~~~~~~~~~~~~~~~
\end{eqnarray}
The extended dimension satisfies a recursion similar to (\ref{dpascal}), from which we find
\eq\label{dimbad}
d^{\nu,\varkappa}_{n,k}=
\left[\!\! \begin{array}{c} \nu-n+1\\ \varkappa-k+1 \end{array}\!\!\right]_{k\alpha-n+1}
{\rm~~~for~~}\alpha > -\infty\,.
\en
A similar formula for $\alpha = -\infty$ requires a further generalisation of Stirling numbers
as in \cite{tsylova}.
To stress dependence on the parameter $\alpha$,
we shall denote the generalised Stirling triangle 
by $G_\alpha$ and denote by ${\cal V}_\alpha$ the set of nonnegative solutions to (\ref{rec1}).
\par By (\ref{dimbad}),
identifying the boundary of $G_\alpha$ is equivalent to finding the limiting regimes
for $\varkappa_\nu$ which entail convergence of certain ratios of the generalised Stirling numbers.
By Lemma \ref{reg} this is reduced to the analysis of possible limits for $V_{n,1}^{\nu,\varkappa}=
d_{n,1}^{\nu,\varkappa}/d^{\nu,\varkappa}$. 
This line seems difficult to pursue, because it requires 
asymptotics of Stirling numbers of distinct types.
Still, there is a much better
formula which 
involves Stirling numbers of a single type:
\begin{lemma} For $\alpha\neq-\infty$ and $\nu\geq n\geq 1,\,\,\nu\geq\varkappa\geq 1$
\begin{equation}\label{vn1}
\V_{n,1}^{\nu,\varkappa}\,d^{n,1}= 
\left(1\bigg/\left[\!\! \begin{array}{c} \nu\\ \varkappa \end{array}\!\!\right]_{\alpha}\right)
\sum_{j=n}^{\nu-\varkappa+1} {\nu-n\choose j-n}
\left[\!\! \begin{array}{c} \nu-j\\ \varkappa-1 \end{array}\!\!\right]_\alpha
(1-\alpha)_{j-1\uparrow}
\end{equation}
where $d^{n,1}=\left[\!\! \begin{array}{c} n\\ 1 \end{array}\!\!\right]_{\alpha}=
(1-\alpha)_{n-1}$.
\end{lemma}
\noindent
{\it Proof.} The left side is the conditional probability 
of $K_n=1$ given $K_\nu=\varkappa$, which is common for all $\V\in {\cal V}_\alpha$. On the other hand,
for partition of the Gibbs form the probability that $\Pi_\nu$ has $\varkappa$ blocks and
the set $[n]$ falls completely in one of the blocks is 
$$\V_{\nu,\varkappa}\sum_{j=n}^{\nu-\varkappa+1} {\nu-n\choose j-n} \W_j \,B_{\nu-j,\varkappa-1}(\W),$$
and to obtain the conditional probability we should divide this by the 
probability $\V_{\nu,\varkappa}B_{\nu,\varkappa}(\W)$ for $\varkappa$ blocks.
Specialising the weights we arrive at (\ref{vn1}).
\endpf

The extension of (\ref{vn1}) to the case $\alpha=-\infty$ is obvious.
The identification of the extremes of ${\cal V}_\alpha$ breaks naturally into cases.

\subsection{Case $-\infty< \alpha<0$}

We claim that the conditions of Lemma \ref{disc} hold with $V(m)=\phi(\alpha, m|\alpha|)$. 
The lemma requires that
\begin{equation}\label{gt1}
\V_{n,m}(m)d^{n,m}={|\alpha|^m m!\over (m|\alpha|)_{n\uparrow}}
\left[\!\! \begin{array}{c} n\\ m \end{array}\!\!\right]_{\alpha}
\end{equation}
goes to $1$ as $n\to\infty$
which seems difficult to check directly. Only in the case $\alpha=-1$ 
this is straightforward due to the 
handy formula for Lah numbers
$$\left[\!\! \begin{array}{c} n\\ m \end{array}\!\!\right]_{-1}={n-1\choose m-1}{n!\over m!}\,.$$

\noindent
Still, (\ref{pvm}) follows trivially from Kingman's representation.
In this case the set $Z$ divides $[0,1]$ into $m$ intervals of random sizes distributed according to the 
symmetric
Dirichlet density 
proportional to $(\xi_1\cdots \xi_m)^{|\alpha|-1}$ on the simplex $\Sigma_j \xi_j=1$.
We have also $\V(m)\to \V(\infty)=\phi(\alpha,\infty\alpha)$ with $\V_{n,k}(\infty)=1(n=k)$.
Thus by Lemma \ref{disc}, the set of extremes is
$\{\phi(\alpha, m|\alpha|), m=1,2,\ldots,\infty\}$. 
As a by-product we have shown that the right side of (\ref{gt1}) approaches $1$ as $n\to\infty$.
\par Another consequence is the following analogue of the Hausdorff moments problem. To
interpret $\V_{n,k}$ as a generalised $(k-1)$th-order difference of the sequence $(\V_{n,1})$,
consider the difference operator $\Delta_\alpha$ which tranforms a sequence $(u_n)$ 
into another sequence
$$(\Delta_{\alpha} u)_n=u_n-(n+1-\alpha)u_{n+1}\,.~~~~~$$
For 
$\V=(\V_{n,k})$ solving (\ref{rec}), setting $(u_n)=(\V_{n,1})$
we have 
\begin{equation}\label{iterate}
\V_{n,k}= (\Delta_{(k-1)\alpha }(\ldots(\Delta_{2 \alpha}(\Delta_{\alpha} u)\ldots)_{n-k+1}\,.
\end{equation}
Note that except for $\alpha = 0$ the operators $\Delta_{j \alpha}$ for different $j$ do not commute.

\begin{corollary} 
\label{moments}
Choose $\alpha<0$.
A sequence 
$(u_n)$ with $u_1=1$ can be represented as
$$
u_n=\sum_{m=1}^\infty {q_m\over (m|\alpha|+1)_{n-1\uparrow}}
$$
for some probability distribution $q$ on $\{1,2,\ldots,\infty\}$ if and only if
the array $(\V_{n,k})$ computed by {\rm (\ref{iterate})} with $(u_n)=(\V_{n,1})$
is nonnegative.
\end{corollary}

\par Another consequence is the asymptotics which does not seem obvious analytically:

$$
\lim_{\nu\to\infty}
{
\left[\!\! \begin{array}{c} \nu-n+1\\ m-k+1 \end{array}\!\!\right]_{k\alpha-n+1}
\over
\left[\!\! \begin{array}{c} \nu\\ m \end{array}\!\!\right]_\alpha~~~~~~
}
= {|\alpha|^k (m)_{k\downarrow}\over (|\alpha|m)_{n\uparrow}}
\qquad {\rm for~~}k=1,\ldots,m\,.
$$

\subsection{Case $\alpha=-\infty$}

This is the limiting case for $\alpha\to-\infty$. Lemma \ref{disc} holds with solutions
$$\V_{n,k}(m) := \phi_{n,k} (-\infty, \infty m) :=
{ (m)_{k\downarrow}\over m^{n}}
\qquad {\rm for~~}k=1,\ldots,m\,.$$
The conditions are checked as in the previous case,
using the fact that $\phi(-\infty, \infty m)$
corresponds to the elementary coupon-collecting partition derived from
a sequence of independent random variables with uniform distribution
on $m$ possible values.

\subsection{Case $\alpha=0$}
That the extremes correspond to the Ewens family can be seen from Lemma \ref{cont}. Take $c_n=\log n$
and recall the well known law of large numbers:  
that for $(0,\theta)$ partition $K_n\sim \theta \log n$ a.s. \cite{abt}.
The fact follows from the representation of $K_n$ as a sum $\xi_1+\ldots +\xi_n$ of independent 
Bernoulli random variables $\xi_j$ with success probability $(\theta+j-1)^{-1}$.


\par We want to see 
how the Ewens family emerges
from the asymptotics of Stirling numbers,
making sense of the statement 
`ESF$(\theta)$ conditioned on $K_n\sim s \log n$ is ESF$(s)$'.
For the Stirling numbers of the first kind there is an asymptotic formula
\cite{abt}
$$
\left[\!\! \begin{array}{c} \nu\\ \varkappa \end{array}\!\!\right]_{0}\sim {\Gamma(\nu)\over\Gamma(\varkappa)\Gamma(1+s)}\,
(\log \nu)^{\varkappa-1}
$$
valid for $\varkappa\sim s\log \nu$, as $\nu\to\infty$
uniformly in $s$ bounded away from $0$ and $\infty$.
Assuming this regime 
for $\varkappa$
some calculus 
shows that (\ref{vn1}) is asymptotic to a Riemann sum for the integral
$$
s\,\int_0^1 x^{n-1} (1-x)^{s-1}{\rm d}x= {(n-1)!\over (s+1)_{n-1\uparrow}}\,.
$$
It is seen that a path is regular if and only if $\varkappa_\nu/\log \nu\to s$ for some $s\in [0,\infty]$
and that the solution induced by such a path has $\V_{n,1}(s)=\phi_{n,1}(0,s)$, hence $\V(s)=\phi(0,s)$.
This identifies the boundary of $G_0$, but 
it is not clear, by this approach, how 
to show that all solutions $\phi(0,\theta)$
are extreme
(which follows from the above law of large numbers).

\par Thus every Gibbs partition with $\W_j=(j-1)!$ is a unique mixture of Ewens' partitions:
$$
p=\int_0^\infty p_{0,\theta} Q({\rm d}\theta)\,,~~~~~~~p\in {\cal P}_0\,
$$
for some probability distribution $Q$ on $[0,\infty]$.
This result was conjectured in \cite{kerov.af, kerov.diss}, and also stated 
without proof in \cite{ki80}, attributed to Frank Kelly.
The analogue of
Corollary \ref{moments} holds with kernel $\theta^k/(\theta)_{n\uparrow}$.

\subsection{Case $0<\alpha<1$}

This is the most interesting case. Lemma \ref{cont} is applied in this case with 
$V(s)=\psi(\alpha|s),\,s\in [0,\infty],$
the law for $(K_n)$ derived from an $(\alpha,\theta)$ partition conditioned on
$K_n\sim s n^\alpha$. 
This is the partition derived by sampling from a Poisson-Kingman random discrete distribution 
denoted in \cite[\S 5.3]{jp.pk} by PK$(\rho_\alpha|t)$ for $t=s^{-\alpha}$.
Here we pursue the connection with Stirling asymptotics.

\par
Let $g_\alpha$ be the density of the Mittag-Leffler distribution, which is determined by the moments
$$
\int_0^\infty x^\beta g_\alpha(x){\rm d}x={\Gamma(\beta+1)\over \Gamma(\beta\alpha+1)}\,,~~~~~\beta>-1\,.
$$
For $k\sim s n^\alpha$, $n\to\infty$ there is the asymptotic formula
$$
\left[\!\! \begin{array}{c} n \\ k \end{array}\!\!\right]_{\alpha}\sim
{\Gamma(n)\over \Gamma(k)} \,n^{-\alpha}\,\alpha^{1-k}\,g_\alpha(s)
$$
which holds uniformly is $s$ bounded away from $0$ and $\infty$.
The formula was derived in
\cite[Section 5]{jp.bm} from a local limit theorem for the stable density 
(note that these Stirling numbers are those of \cite{jp.bm} multiplied by $\alpha^{-k}$).
In \cite[Theorem 11]{fl-maps} the 
formula was concluded by the singularity
analysis of the generating function 
(\ref{genfun}), 
and in \cite[Corollary 12]{fl-urns} the formula appeared in connection
with an urn model similar to that in \cite{tsylova}.

\par Substituting this approximation into (\ref{vn1}), using
$\Gamma(\nu+\beta)/\Gamma(\nu)\sim \nu^{\beta}$ and introducing the variable $y=j/\nu$ we arrive 
at a Riemann sum for
\begin{equation}\label{lim-v}
\V_{n,1}(s) \,(1-\alpha)_{n-1\uparrow}= {s\alpha\over \Gamma(1-\alpha) g_\alpha(s)}\int_0^1 y^{n-1-\alpha}
(1-y)^{-1-\alpha}
g_\alpha(s(1-y)^{-\alpha})\,{\rm d}y.
\end{equation}
Using the change of variable $t=s^{-1/\alpha}$ and the formula
$$
f_\alpha(y)=\alpha y^{-1-\alpha}\,g_\alpha (y^{-\alpha})
$$
connecting $g_\alpha$ to the stable density $f_\alpha$, whose Laplace transform at $\lambda$ is $\exp(-\lambda^\alpha)$,
we get 
$$
\V_{n,1}(s) \,(1-\alpha)_{n-1\uparrow}=
{\alpha\over \Gamma(1-\alpha)t^\alpha f_\alpha(t)} \int_0^1 y^{n-1-\alpha} f_\alpha(t(1-y)){\rm d}y\,,
$$
which is an instance of \cite[Equation (66)]{jp.pk}.
It follows that a path is regular if and only if $\varkappa_\nu/\nu^\alpha$ converges,
in which case the induced solution is $\psi(\alpha|s)$.
This identifies the corresponding partition
as a Poisson-Kingman partition, as indicated above. 

\par Thus the Stirling asymptotics enable to determine the boundary of $G_\alpha$, but we do not see how they 
imply that all boundary solutions are extreme. 
But this is obvious from the law of large numbers: because 
for $(\alpha,0)$-partition $K_n/n^\alpha$ has a strong limit, the same is true
for the $(\alpha|s)$ partition, obtained by conditioning. 
The paintbox $Z$ for the
$(\alpha|s)$ partition
can be represented by the division of $[0,1]$ into excursions of a Bessel 
bridge (or Bessel process) conditioned on its local time at zero being equal to $s$.

\vskip0.5cm

\par We summarise conclusions of this section in the following theorem, which is our main result.

\begin{theorem} 
\label{mainth}
Each Gibbs partition of fixed type $\alpha\in [-\infty,1[\,$ is a unique probability mixture of the
extreme partitions of this type, which are
\begin{enumerate}
\item[{\rm (i)}] $(\alpha,|\alpha| m)$-partitions with $m=0,1,\ldots,\infty$, for $\alpha\in \,[-\infty,0[\,;$
\item[{\rm (ii)}] the Ewens $(0,\theta)$-partitions with $\theta\in [0,\infty]$, for $\alpha=0$\,;
\item[{\rm (iii)}] the Poisson-Kingman $(\alpha|s)$-partitions with $s\in [0,\infty]$, for $\alpha\in \,]0,1[\,$.
\end{enumerate}
\end{theorem}

\vskip0.5cm
\noindent
{\bf Acknowledgement} We are indebted for Philippe Flajolet and Grigori Olshanski 
for alerting us of the connections with their work.

\def\cprime{$'$} \def\polhk#1{\setbox0=\hbox{#1}{\ooalign{\hidewidth
\lower1.5ex\hbox{`}\hidewidth\crcr\unhbox0}}} \def\cprime{$'$}
\def\cprime{$'$} \def\cprime{$'$}
\def\polhk#1{\setbox0=\hbox{#1}{\ooalign{\hidewidth
\lower1.5ex\hbox{`}\hidewidth\crcr\unhbox0}}} \def\cprime{$'$}
\def\cprime{$'$} \def\polhk#1{\setbox0=\hbox{#1}{\ooalign{\hidewidth
\lower1.5ex\hbox{`}\hidewidth\crcr\unhbox0}}} \def\cprime{$'$}
\def\cprime{$'$} \def\cydot{\leavevmode\raise.4ex\hbox{.}} \def\cprime{$'$}
\def\cprime{$'$} \def\cprime{$'$} \def\cprime{$'$}

\end{document}